\documentclass[reqno]{amsart}
\usepackage{amsmath}
\usepackage{amscd}
\usepackage{amssymb}
\usepackage{graphics}
\usepackage{latexsym}
\usepackage{hyperref}

\theoremstyle{plain}
\newtheorem{theorem}{Theorem}[section]
\newtheorem{thm}[theorem]{Theorem}
\newtheorem{cor}[theorem]{Corollary}
\newtheorem{lem}[theorem]{Lemma}
\newtheorem{prop}[theorem]{Proposition}

\theoremstyle{definition}
\newtheorem{defn}[theorem]{Definition}
\newtheorem{ques}[theorem]{Question}

\theoremstyle{remark}

\newcommand{\ZZ}{\mathbb{Z}}
\newcommand{\CC}{\mathbb{C}}

\newcommand{\AAA}{\mathbb{A}}
\newcommand{\PP}{\mathbb{P}}

\newcommand{\SP}{\text{Spec }}

\newcommand{\marpar}[1]{}

\newcommand{\lt}{\left}
\newcommand{\rt}{\right}

\newsavebox{\sembox}
\newlength{\semwidth}
\newlength{\boxwidth}

\newsavebox{\semrbox}
\newlength{\semrwidth}
\newlength{\boxrwidth}

\newcommand{\mc}{\mathcal}
\newcommand{\mb}{\mathbf}

\newcommand{\OO}{\mathcal O}

\def\PP{{\mathbb P}}

\begin{document}

\title[A pencil of index one]
{A pencil of Enriques surfaces of index one with no section}

\author[Starr]{Jason Michael Starr}
\address{Department of Mathematics \\
  Massachusetts Institute of Technology \\ Cambridge MA 02139}
\email{jstarr@math.mit.edu} 
\date{\today}

\begin{abstract}
  Monodromy arguments and deformation-and-specialization are used to
  prove existence of a pencil of Enriques surfaces with no section and
  \emph{index 1}.  The same technique completes the strategy
  from ~\cite[\S 7.3]{GHMS} proving the family of \emph{witness curves
    for dimension $d$} depends on the integer $d$.
\end{abstract}

\maketitle

\section{Introduction} ~\label{intro}
% \marpar{intro}

\noindent
This paper uses monodromy and deformation-and-specialization to answer
some questions related to ~\cite{GHMS}.  Theorem~\ref{thm-main} gives
a new, elementary proof of existence of a pencil of Enriques surfaces
over $\CC$ with no section, which moreover has \emph{index 1}.
Proposition~\ref{prop-sec7} completes the strategy from ~\cite[\S
7.3]{GHMS} proving the family $\mc{H}_d$ of \emph{witness curves}
depends on the relative dimension $d$.

\medskip\noindent
The main theorem of ~\cite{GHS} proves a rationally connected variety
defined over the function field of a curve over a characteristic $0$
algebraically closed field has a rational point.  A converse 
is proved in ~\cite{GHMS}; in particular ~\cite[Cor. 1.4]{GHMS} proves
there is an Enriques surface without a rational point that is defined
over the function field of a curve (answering a question of Serre
~\cite[p. 153]{GtoS}).  Subsequently Lafon~\cite{Lafon} gave an
\emph{explicit} pencil of Enriques surfaces defined over $\ZZ[1/2]$
whose base-change to any field of characteristic $\neq 2$ has no
rational point.  H\'el\`ene Esnault asked about the index of
Enriques surfaces without a rational point.

\medskip
\begin{defn} \label{defn-index}
% \marpar{defn-index}
  Let $X$ be a finite type scheme, algebraic space, algebraic stack,
  etc. over a field $K$.  The \emph{index} and the \emph{minimal
    degree} are,
$$
\begin{array}{ccc}
I(K,X) & = & \text{gcd}\{ [L:K] | X(L) \neq \emptyset \}, \\
M(K,X) & = & \min \{ [L:K] | X(L) \neq \emptyset \}.
\end{array}
$$
\end{defn}

\medskip\noindent
H\'el\`ene Esnault asked, essentially, what is the possible index of
an Enriques surface defined over a function field of a curve.  In
Lafon's example, $M(K,X_{K})=I(K,X_{K})=2$.  In ~\cite{GHMS} the index
is not computed, but likely there also $I(K,X_K) > 1$.

\begin{ques}[Esnault] \label{ques-Esnault}
% \marpar{ques-Esnault}
  If $X$ is an Enriques surface defined over a function field of a
  curve $K$ with no $K$-point, is $I(K,X) > 1$?
\end{ques}

\medskip\noindent
This has to do with whether there is an obstruction to $K$-points in
Galois cohomology.  If so and if the obstruction is compatible with
restriction and corestriction, the order of the obstruction divides
$I(K,X)$.  So if there is a cohomological obstruction ``explaining''
non-existence of $K$-points, then $I(K,X_K)>1$.  The main result
proves there is an Enriques surface with no $K$-point whose index is
$1$.

\begin{thm} \label{thm-main}
% \marpar{thm-main}
  Let $k$ be an algebraically closed field with $\text{char}(k)\neq 2,
  3$ that is ``sufficiently big'', e.g. uncountable.  There
  exists a flat, projective $k$-morphism
  $\pi:\mc{X} \rightarrow \PP^1_k$ with the following properties,
\begin{enumerate}
\item[(i)] 
  the geometric generic fiber of $\pi$ is a smooth Enriques surface,
\item[(ii)]    
  the invertible sheaf $\pi_*[\omega_\pi^{\otimes 2}]$ has degree $6$,
\item[(iii)]  
  for the function field $K$ of $\PP^1_k$ and the generic fiber $X_K$
  of $\pi$, $I(K,X_K)=1$ and $M(K,X_K)=3$.
\end{enumerate} 
Moreover every ``very general'' Enriques surface over $k$ is a fiber
of such a family.
\end{thm}

\medskip\noindent
The method is simple.  Over $\PP^1$ a family of surfaces is given
whose monodromy group acts as the full group of symmetries of the dual
graph of the geometric generic fiber -- which is the 2-skeleton of a
cube.  There is an action of $\ZZ/2\ZZ$ acting fiberwise, and the
quotient is a pencil $\mc{X}/\PP^1$ of ``Enriques surfaces''.  The $8$
vertices of the cube give a degree $4$ multi-section of the pencil.
The $6$ faces of the cube give a degree $3$ multi-section of the
pencil.  By monodromy considerations every multi-section of $X$ has
degree $\geq 3$.  The pencil $X$ together with the degree $3$ and
degree $4$ multi-sections deforms to a pencil whose geometric generic
fiber is a smooth Enriques surface.  For a general such deformation,
$M(K,X_K)=3$ and $I(K,X_K)=1$.

\medskip\noindent 
The same method gives pencils of degree $d$ hypersurfaces with minimal
degree $d$, which is used to complete the argument from ~\cite[Section
7.3]{GHMS}.

\begin{prop} \label{prop-sec7}
% \marpar{prop-sec7}  
  Let $B$ be a normal, projective variety of dimension $\geq 2$ and
  let $M$ be an irreducible family of irreducible curves dominating
  $B$ (i.e., the morphism from the total space of the family of curves
  to $B$ is dominant).  There is an integer $d$ such that $M$ is not a
  \emph{witness family for dimension $d$}, i.e., there is a
  projective, dominant morphism of relative dimension $d$,
  $\pi:\mc{X}\rightarrow B$, whose restriction to each curve of $M$
  has a section, but whose restriction to some smooth curve in $B$ has
  no section.
\end{prop}

\medskip\noindent
\textbf{Acknowledgments:}  
I am grateful to Igor Dolgachev, H\'el\`ene Esnault, 
Mike Stillman and Harry Tamvakis for
useful conversations.  This paper originated in meetings with Tom
Graber, Joe Harris and Barry Mazur to whom I am very grateful.  I was
partially supported by NSF Grant DMS-$0201423$.

\medskip\noindent
\section{The construction for hypersurfaces} \label{sec-constr} 
% \marpar{sec-constr} 

\noindent
Let $d,n >0$ be integers, let $k$ be a field, and let $V$ be a
$k$-vector space of dimension $n+1$.  Degree $d$ hypersurfaces in
$\PP(V)$ are parametrized by the projective space,
$$
\PP\text{Sym}^d(V^\vee) = \text{Proj} \bigoplus_i
\text{Sym}^i(\text{Syt}^d(V)), 
$$ 
where $\text{Syt}^d(V)$ is the vector space of symmetric tensors in
$\otimes^d V$.

\medskip\noindent
Let $B,C$ be $k$-curves isomorphic to $\PP^1_k$.  There exists a
degree $d$, separably-generated $k$-morphism $f:C\rightarrow B$ such
that $\text{Gal}(k(C)/k(B))$ is the full symmetric group
$\mathfrak{S}_d$.  This is straightforward in every characteristic --
in characteristic $0$ any morphism with simple branching will do.

\medskip\noindent
Let $g:C\rightarrow \PP(V^\vee)$ be a closed immersion whose image is
a rational normal curve of degree $n$.  Consider the pullback of the
tautological surjection, $V \otimes_k \OO_C \rightarrow g^*\OO(1)$.
By adjointness, there is a map $\beta: V\otimes_k \OO_B \rightarrow
f_*(g^*\OO(1))$.  For every locally free $\OO_C$-module $\mc{E}$ there
is the \emph{norm sheaf} on $B$,
\begin{equation*}
\text{Nm}_f(\mc{E}) =
\textit{Hom}_{\OO_B}(\bigwedge^d(f_*\OO_C), \bigwedge^d(f_*\mc{E}) ),
\end{equation*}
together with the \emph{norm map} of 
$\OO_B$-modules,
\begin{equation*}
\alpha'_{\mc{E}}: \bigotimes^d (f_*\mc{E}) \rightarrow
\text{Nm}_f(\mc{E}), \ \ e_1\otimes \dots \otimes e_d \mapsto \lt(
c_1\wedge \dots \wedge c_d 
\mapsto (c_1\cdot e_1) \wedge \dots \wedge (c_d\cdot e_d) \rt), 
\end{equation*}
for $e_1\otimes\dots \otimes e_d\in \bigotimes^d(f_*\mc{E})$ and
$c_1\wedge \dots \wedge c_d \in \bigwedge^d(f_*\OO_B)$.  Only the
restriction to the subsheaf of symmetric tensors is needed,
$\alpha_{\mc{E}}: \text{Syt}^d(f_*\mc{E}) \rightarrow
\text{Nm}_f(\mc{E})$.  In particular, $\text{Nm}_f(\OO_C) = \OO_B$ and
$\alpha_{\OO_C}(b\otimes\dots \otimes b) \in \OO_B$ is the usual norm
of $b\in f_*\OO_C$.

\medskip\noindent
Denote by $\gamma$ the composition,
$$
\text{Syt}^d(V)\otimes_k \OO_B \xrightarrow{\text{Syt}^d(\beta)}
\text{Syt}^d(f_*g^*\OO(1)) \xrightarrow{\alpha_{g^*\OO(1)}}
  \text{Nm}_f(g^* \OO(1)).
$$
Because $\beta$ is surjective, also $\gamma$ is surjective.  So there
is an induced morphism $h:B \rightarrow \PP \text{Sym}^d(V^\vee)$.  For
every geometric point $b\in B$ whose fiber $f^{-1}(b)$ is a reduced
set $\{c_1,\dots,c_d\}$, $h(b) = [g(c_1)\times \cdots \times g(c_d)]$.
The degree of $\text{Nm}_f(g^*\OO(1))$, and thus the degree of $h$, is
$n$

\medskip\noindent
Denote by $\mc{X}_h\subset B\times \PP(V)$ the preimage under
$(h,\text{Id})$ of the universal hypersurface in $\PP\text{Sym}^d
(V^\vee)\times \PP(V)$, and by $\pi:\mc{X}_h \rightarrow B$ the
projection.  Let $m=\min(d,n)$ and let $S_{d,n} \subset \ZZ_{\geq 0}$
denote the additive semigroup generated by $\binom{d}{i}$ for
$i=1,\dots,m$.  Denote $K=k(B)$ and denote by $\mc{X}_{h,K}$ the
generic fiber of $\pi$.

\begin{prop} \label{prop-Sdn}
% \marpar{prop-Sdn}  
  Every irreducible multi-section of $\pi$ has degree divisible by
  $\binom{d}{i}$ for $i=1,\dots,m$.  The degree of every multi-section
  is in $S_{d,n}$.  In particular, if $d > n$ then $M(K,\mc{X}_{h,K})
  = d$ and $I(K,\mc{X}_{h,K})$ is divisible by
  $\text{gcd}(d,\binom{d}{2},\dots,\binom{d}{n})$.
\end{prop}  

\begin{proof}  
  Denote by $U \subset B$ the largest open subset over which $f$ is
  \'etale and define $W=f^{-1}(U)$.  For each $i=1,\dots,m$, denote by
  $W_i/U$ the relative Hilbert scheme $\text{Hilb}^i_{W/U}$.  Because
  $W$ is \'etale over $U$, the fiber of $f$ over a geometric point
  $b$ of $B$ is a set of $d$ distinct points, $f^{-1}(b)=
  \{c_1,\dots,c_d\}$, and 
  the fiber of $\text{Hilb}^i_{W/U}$ is the set of subsets of
  $f^{-1}(b)$ of size $i$. 
Every
  geometric fiber of $\mc{X}_h\times_B U \rightarrow U$ is 
union of $d$
  hyperplanes.  Denote by,
\begin{equation*}
\mc{X}_h \times_B U = \mc{X}_h^{1} \sqcup \mc{X}_h^{2} \sqcup \dots \sqcup
\mc{X}_h^{n},
\end{equation*}
the locally closed stratification where $\mc{X}_h^i$ is the set of
points $x$ in precisely $i$ irreducible components of the geometric
fiber $\mc{X}_h \otimes_{\OO_B} \overline{\kappa}(\pi(x))$.  Because
every finite subset of distinct closed points on a rational normal curve over
an algebraically closed field is
in linearly general position, $\mc{X}_h^i = \emptyset$ for $i > m$;
in particular every geometric fiber of $\mc{X}_h\times_B U \rightarrow
U$ is a simple normal crossings variety.  
For each
$i=1,\dots,m$ the morphism $\mc{X}^i_h \rightarrow U$ factors as an
$\AAA^{n-i}$-bundle over $W_i$ over $U$.  
The generic point of every irreducible
multi-section is contained in $\mc{X}^i_h$ for some $i=1,\dots,m$.
Because $\text{Gal}(k(C)/k(B))$ is $\mathfrak{S}_d$, $W_i$ is
irreducible.  Therefore the degree of the multi-section is divisible
by $\text{deg}(k(W_i)/k(U)) = \binom{d}{i}$.  So the degree of every
multi-section, irreducible or not, is in $S_{d,n}$.  Moreover, the
intersection of $\mc{X}_{h,K}$ with a general line in $\PP(V\otimes_k
K)$ is a degree $d$ multi-section, so $M(K,\mc{X}_{h,K})=d$.  
\end{proof}

\medskip\noindent
Let $H_n \subset \text{Hom}(B,\PP \text{Sym}^d(V^\vee))$ denote the
irreducible component of morphisms of degree $n$.  Denote by
$\mc{X}\rightarrow H_n\times B$ the pullback by the universal morphism
of the universal hypersurface in $\PP \text{Sym}^d(V^\vee)\times
\PP(V)$.  
For every field
$k'$ and every $[h]\in H_n(k')$, denote by $\mc{X}_h$ the restriction
of $\mc{X}$ to $\SP(k') \times B$,
by $K'$ the function field
$k'(B)$, and by $\mc{X}_{h,K'}$ the generic fiber of the projection to
$B$.

\begin{cor} \label{cor-Sdn}
% \marpar{cor-Sdn}  
  Assume $d > n$.  In $H_n$ there is a countable intersection of open
  dense subsets such that for every $[h]$ in this set,
  $M(K',\mc{X}_{h,K'})=d$ and $I(K',\mc{X}_{h,K'})$ is divisible by
  $\text{gcd}(d,\dots,\binom{d}{n})$. In particular this holds for the
  geometric generic point of $H_n$.
\end{cor}

\begin{proof}    
  The subset $H_n^{\text{good}}\subset H_n$ where $M(K',\mc{X}_{K'})
  \geq d$ and $\text{gcd}(d,\dots,\binom{d}{n}) \mid
  I(K',\mc{X}_{h,K'})$ is a countable intersection of open subsets by
  standard Hilbert scheme arguments: the complement of this set is the
  union over the countably many Hilbert polynomials $P(t)$ of
  multi-sections of degree $<d$ or not divisible by
  $\text{gcd}(d,\dots,\binom{d}{n})$ of the closed image in $H_n$ of
  the relative Hilbert scheme $\text{Hilb}^{P(t)}_{\mc{X}/H_n}$.  By
  Proposition~\ref{prop-Sdn} $H_n^{\text{good}}$ is nonempty,
  therefore it is a countable intersection of open \emph{dense}
  subsets.  Of course the intersection of $\mc{X}_{h,K'}$ with a
  general line in $\PP(V\otimes_k K')$ gives a multi-section of degree
  $d$, therefore $H_n^{\text{good}}$ is actually the set where
  $M(K',\mc{X}_{K'})=d$ and $\text{gcd}(d,\dots,\binom{d}{n}) \mid
  I(K',\mc{X}_{h,K'})$.
\end{proof}

\subsection{Proof of Proposition~\ref{prop-sec7}} \label{subsec-proof}
% \marpar{subsec-proof}
Let $k$ be an uncountable, algebraically closed field.  The main case
of Proposition~\ref{prop-sec7} is $B=\PP^1_k\times \PP^1_k$ and $M$ is
the complete linear system $|\OO(a,b)|$.  Assume first that one of
$a,b=0$, say $b=0$.  Let $f:Y\rightarrow \PP^1_k$ be a finite,
separably-generated morphism of irreducible curves of degree $>1$, and
let $\mc{X} = Y\times \PP^1$ with projection $\pi=(f,\text{Id})$.
Every divisor in $|\OO(a,0)|$ is a union of fibers of $\text{pr}_1$,
so the restriction of $\pi$ has a section.  The restriction of $\pi$
over every fiber of $\text{pr}_2$ is just $f$, and so has no rational
section.  Thus assume $a,b>0$.

\medskip\noindent
Define $n= 4ab$ and $d=n-1$.  Let $V$ be a $k$-vector space of
dimension $n+1$.  Let $C \subset \PP^1\times\PP^1$ 
be a smooth curve in the linear system
$|\OO(1,2b)|$.  By Corollary~\ref{cor-Sdn}, there exists a
closed immersion of degree $n$, $h:C\rightarrow
\PP\text{Sym}^d(V^\vee)$, such that $M(k(C),\mc{X}_{h,k(C)}) = d>1$.
Of course $h$ extends to a closed immersion $j:\PP^1\times \PP^1
\rightarrow \PP \text{Sym}^d(V^\vee)$ such that $j^*\OO(1) =
\OO(2a-1,2b)$; after all, $H^0(\PP^1\times\PP^1, \OO(2a-1,2b))
\rightarrow H^0(C,\OO_C(n))$ is surjective.  Define $\pi:\mc{X}
\rightarrow \PP^1\times \PP^1$ to be the base-change by $j$ of the universal
family of degree $d$ hypersurfaces in $\PP(V)$.  By construction, the
restriction over $C$ has no section.

\medskip\noindent  
Every divisor in $|\OO(a,b)|$ is a curve in $\PP\text{Sym}^d(V^\vee)$
of degree $n-b$ whose span is a linear system of hypersurfaces in
$\PP(V)$ of (projective) dimension $\leq n-b-(a-1)(b-1)$.
Since $n-b<n$, this
linear system has basepoints giving sections of the restriction of
$\mc{X}$ to the divisor.  This proves Proposition~\ref{prop-sec7}
for $B=\PP^1\times \PP^1$ and $M=|\OO(a,b)|$.

\medskip\noindent
Let $B$ be a normal, projective variety of dimension $\geq 2$ and let
$M$ be an irreducible family of irreducible curves dominating $B$.
There exists a smooth open subset $U\subset B$ whose complement has
codimension $\geq 2$ and a dominant morphism $g:U\rightarrow
\PP^1\times \PP^1$.  Intersecting $U$ with general hyperplanes, there
exists an irreducible closed subset $Z\subset U$ such that
$g|_Z:Z\rightarrow \PP^1\times \PP^1$ is generically finite of some
degree $e>0$.  For the geometric generic point of $M$, the intersection
of the corresponding curve with $U$ is nonempty, and the closure of
the image under $f$ is a divisor in the linear system $|\OO(a',b')|$
for some integers $a',b'$.  Let $a\geq a'$, and  $b\geq b'$ be integers such
that $4ab > e+1$.  There exists a projective, dominant morphism
$\pi:\mc{X} \rightarrow \PP^1\times \PP^1$ whose restriction over
every divisor in $|\OO(a,b)|$ has a section, but whose restriction
over a general divisor in $|\OO(1,2b)|$ has minimal degree $4ab-1$.

\medskip\noindent
Define $\mc{X}_B\subset B\times \mc{X}$ to be the closure of
$U\times_{\PP^1 \times \PP^1} \mc{X}$.  Then
$\pi_B:\mc{X}_B\rightarrow B$ is a projective dominant morphism.  For
the geometric generic point of $M$, the restriction of $\pi_B$ to the
curve has a section because the restriction of $\pi$ to the image in
$\PP^1\times \PP^1$ has a section.  Let $C_B\subset Z$ be the preimage
of a general curve $C$ in $|\OO(1,2b)|$.  The morphism $C_B\rightarrow
C$ has degree $e < 4ab-1$.  Because every multi-section of $\pi$
over $C$ has degree $\geq 4ab-1$, $\pi_B$ has no section over
$C_B$.

\section{The construction for Enriques surfaces} \label{sec-conES}
% \marpar{sec-conES}

\noindent
Let $k$ be a field of characteristic $\neq 2,3$, and let $V_+$ and
$V_-$ be $3$-dimensional $k$-vector spaces.  Denote $V=V_+\oplus V_-$
and $V' = \text{Sym}^2(V^\vee_+)\oplus \text{Sym}^2(V^\vee_-)$.
Denote $G=\text{Grass}(3,V')$, parametrizing $3$-dimensional
\emph{subspaces} of $V'$.  This is a parameter space for Enriques
surfaces.  There are 2 descriptions of the universal family, each
useful.  First, let $\pi_Z:Z\rightarrow \PP(V_+)\times \PP(V_-)$ be
the projective bundle of the locally free sheaf
$\text{pr}_+^*\OO_{\PP(V_+)}(-2)\oplus
\text{pr}_-^*\OO_{\PP(V_-)}(-2)$.  A general complete intersection of
3 divisors in $|\OO_Z(1)|$ is an Enriques surface.  Because
$H^0(Z,\OO_Z(1)) = V'$, the parameter space for these complete
intersections is $G$.  Second, $G$ parametrizes complete intersections
in $\PP(V)$ of 3 quadric divisors that are invariant under the
involution $\iota$ of $\PP(V)$ whose $(-1)$-eigenspace is $V_-$ and
whose $(+1)$-eigenspace is $V_+$.  A general such complete
intersection is a K3 surface on which $\iota$ acts as a
fixed-point-free involution; the quotient by $\iota$ is an Enriques
surface.  The two descriptions are equivalent: the involution extends
to an involution $\widetilde{\iota}$ on the blowing up
$\widetilde{\PP(V)}$ of $\PP(V)$ along $\PP(V_+)\cup \PP(V_-)$ and the
quotient is $Z$.  Denote by $\mc{X}\rightarrow G$ the universal family
of Enriques surfaces, and denote by $\mc{Y}\rightarrow G$ the
universal family of K3 covers.

\medskip\noindent
Let $B, C, D$ be $k$-curves isomorphic to $\PP^1_k$.  There exists a
degree 2, separably-generated morphism $g:D\rightarrow C$ and a degree
3, separably-generated morphism $f:C\rightarrow B$ such that
$\text{Gal}(k(D)/k(B))$ is the full wreath product
$\mathfrak{W}_{3,2}$, i.e., the semidirect product $(\mathfrak{S}_2)^3
\rtimes \mathfrak{S}_3$.  In characteristic $0$, this holds whenever
$g$ and $f$ have simple branching and the branch points of $g$ are in
distinct, reduced fibers of $f$.  There is an involution $\iota_D$ of
$D$ commuting with $g$.

\medskip\noindent
Let $j:D\rightarrow \PP(V^\vee)$ be a closed immersion equivariant for
$\iota_D$ and $\iota$ whose image is a rational normal curve of degree
$5$.  By the construction in Section~\ref{sec-constr}, there is an
associated morphism $i:C\rightarrow \PP \text{Sym}^2(V^\vee)$.
Because $j$ is equivariant, $i$ factors through $\PP(V')$.  By a
straightforward computation, $i^*\OO(1) = \text{Nm}_g(j^*\OO(1)) \cong
\OO_C(5)$.  The pushforward by $f_*$ of the pullback by $i^*$ of the
tautological surjection is a surjection $(V')^\vee\otimes \OO_B
\rightarrow f_*i^*\OO(1)$.  The sheaf $f_*i^*\OO(1)$ is locally free,
in fact $f_*i^*\OO(1) \cong f_*\OO_C(5) \cong \OO_B(1)^3$, so there is
an induced morphism $h:B\rightarrow G$.  Denote by
$\pi_h:\mc{X}_h\rightarrow B$ and $\rho_h:\mc{Y}_h \rightarrow B$ the
base-change by $h$ of $\mc{X}$ and $\mc{Y}$.  Denote $K=k(B)$ and denote
by $\mc{X}_{h,K}$ the generic fiber of $\pi_h$.

\begin{prop} \label{prop-W32}
% \marpar{prop-W32}  
  Every irreducible multi-section of $\pi_h$ has degree divisible by
  $3$ or $4$.  In
  particular $M(K,\mc{X}_{h,K}) = 3$.
\end{prop}

\begin{proof}  
  Denote by $U\subset B$ the open set over which $f\circ g$ is
  \'etale, and denote by $W\subset D$ the preimage of $U$.  Denote by
  $c:\widetilde{W}\rightarrow U$ the Galois closure of $W/U$.  Then
  $c^* f_*\OO_C|_U \cong
  \OO_{\widetilde{W}}\{\mb{a}_1,\mb{a}_2,\mb{a}_3\}$ for idempotents
  $\mb{a}_p$, $p=1,2,3$.  And $c^* g_*f_*\OO_D|_U \cong
  \OO_{\widetilde{W}}\{
  \mb{b}_{1,1},\mb{b}_{1,2},\mb{b}_{2,1},\mb{b}_{2,2},
  \mb{b}_{3,1},\mb{b}_{3,2}\}$ for idempotents $\mb{b}_{p,q}$,
  $p=1,2,3$, $q=1,2$.  Of course $\mb{a}_p \mapsto
  \mb{b}_{p,1}+\mb{b}_{p,2}$, $p=1,2,3$.  The action of the Galois
  group $\mathfrak{W}_{3,2}$ on $\mb{a}_p$ is by the symmetric group
  $\mathfrak{S}_3$, and the action on $\mb{b}_{p,q}$ is the standard
  representation of the wreath product.

\medskip\noindent
For each $p=1,2,3$ and $q=1,2$, denote by
$j_{p,q}:\widetilde{W}\rightarrow \PP(V^\vee)$ the morphism obtained
by composing the idempotent $\mb{b}_{p,q}:\widetilde{W} \rightarrow
\widetilde{W}\times_U W$ with the basechange of $j$.  In particular,
$\iota\circ j_{p,1} = j_{p,2}$.  Denote by $\Lambda_{p,q} \subset
\widetilde{W}\times \PP(V)$ the pullback by $(j_{p,q},\text{Id})$ of
the universal hyperplane.  Denote by $\mc{Y}_{\widetilde{W}}$ the
base-change to $\widetilde{W}$ of $\mc{Y}_h$.  Then,
$$
\mc{Y}_{\widetilde{W}} = \bigcup_{(q_1,q_2,q_3)\in\{1,2\}^3}
(\Lambda_{1,q_1}\cap \Lambda_{2,q_2}\cap \Lambda_{3,q_3}).
$$

\medskip\noindent
There is a locally closed stratification,
$$
\mc{Y}_{\widetilde{W}} = \mc{Y}_{\widetilde{W}}^3 \sqcup
\mc{Y}_{\widetilde{W}}^4 \sqcup \mc{Y}_{\widetilde{W}}^5,
$$
where $\mc{Y}_{\widetilde{W}}^l$ is the set of points lying in the
intersection of precisely $l$ of the $\Lambda_{p,q}$.  The stratum
$\mc{Y}_{\widetilde{W}}^3$ is the union of 8 connected, open subsets,
$$
\Lambda_{(q_1,q_2,q_3)} \subset
(\Lambda_{1,q_1}\cap \Lambda_{2,q_2}\cap \Lambda_{3,q_3}),
$$
for $q_1,q_2,q_3\in \{1,2\}$.  Each connected component is a dense
open subset of a $\PP^2$-bundle over $\widetilde{W}$.  The stratum
$\mc{Y}_{\widetilde{W}}^4$ is the union of 12 connected, open subsets,
$$
\begin{array}{ccc}
\Lambda_{(*,q_2,q_3)} & \subset & (\Lambda_{1,1}\cap \Lambda_{1,2})
\cap \Lambda_{2,q_2} \cap \Lambda_{3,q_3}, \\
\Lambda_{(q_1,*,q_3)} & \subset & \Lambda_{1,q_1} \cap
(\Lambda_{2,1}\cap \Lambda_{2,2}) \cap \Lambda_{3,q_3}, \\
\Lambda_{(q_1,q_2,*)} & \subset & \Lambda_{1,q_1} \cap \Lambda_{2,q_2}
\cap (\Lambda_{3,1}\cap \Lambda_{3,2})
\end{array}
$$
for $q_1,q_2,q_3\in \{1,2\}$.  Each connected component is a dense
open subset of a $\PP^1$-bundle over $\widetilde{W}$.  Finally
$\mc{Y}_{\widetilde{W}}^5$ is the union of 6 connected sets,
$$
\begin{array}{ccc}
\Lambda_{(*,*,q_3)} & = & (\Lambda_{1,1} \cap \Lambda_{1,2})
\cap (\Lambda_{2,1}\cap \Lambda_{2,2}) \cap \Lambda_{3,q_3}, \\
\Lambda_{(*,q_2,*)} & = & (\Lambda_{1,1} \cap \Lambda_{1,2}) \cap
\Lambda_{2,q_2} \cap (\Lambda_{3,1} \cap \Lambda_{3,2}), \\
\Lambda_{(q_1,*,*)} & = & \Lambda_{1,q_1} \cap (\Lambda_{2,1} \cap
\Lambda_{2,2}) \cap (\Lambda_{3,1} \cap \Lambda_{3,2})
\end{array}
$$
for $q_1,q_2,q_3  \in \{1,2\}$.   Each  connected component  projects
isomorphically to $\widetilde{W}$.

\medskip\noindent
There is a bijection between multi-sections of $\mc{Y}_h$ over $U$ and
Galois invariant multi-sections of $\mc{Y}_{\widetilde{W}}$ over
$\widetilde{W}$.  An irreducible multi-section of $\mc{Y}_h$
determines a multi-section of $\mc{Y}_{\widetilde{W}}$ contained in a
single stratum $\mc{Y}^l_{\widetilde{W}}$.  The action of the Galois
group $\mathfrak{W}_{3,2}$ on the connected components of
$\mc{Y}^l_{\widetilde{W}}$ is the obvious one; in particular, it acts
transitively on the set of connected components.  So every Galois
invariant multi-section in $\mc{Y}^3_{\widetilde{W}}$ has degree
divisible by $8$, every Galois invariant multi-section in
$\mc{Y}^4_{\widetilde{W}}$ has degree divisible by $12$, and every
Galois invariant multi-section in $\mc{Y}^5_{\widetilde{W}}$ has
degree divisible by $6$.  Therefore every irreducible multi-section of
$\mc{Y}_h$ has degree divisible by $8$ or $6$.  Because $\mc{Y}_h$ is
a double-cover of $\mc{X}_h$, every irreducible multi-section of
$\mc{X}_h$ has degree divisible by $4$ or $3$.  In particular, the
minimal degree of a multi-section of $\mc{X}_h$ is $3$.
\end{proof}

\medskip\noindent
Because $f_*i^*\OO(1) \cong \OO_B(1)^3$, the scheme $\mc{X}_h \subset
B \times Z$ is a complete intersection of 3 divisors in the linear
system $|\text{pr}_B^*\OO_B(1) \otimes \text{pr}_Z^* \OO_Z(1)|$.  A
general deformation of this complete intersection is a pencil of
Enriques surfaces satisfying Theorem~\ref{thm-main} (i) and (ii) with
$M(K,X_K) \geq 3$, $I(K,X_K) \mid 4$ (this is valid so long as
$\text{char}(k)\neq 2$).  For (iii), it is necessary to deform the
pencil together with the degree $3$ multi-section.  This requires a
bit more work, and the hypothesis $\text{char}(k) \neq 2,3$.

\medskip\noindent
The stratum $\mc{Y}^5_{\widetilde{W}}$ is Galois invariant and
determines a degree $3$ multi-section of $\mc{X}_h$.  As a
$\mathfrak{W}_{3,2}$-equivariant morphism to $\widetilde{W}$,
$\mc{Y}_{\widetilde{W}}^5$ is just the base-change of $D$, and the
morphism $\mc{Y}_{\widetilde{W}}^5 \rightarrow \PP(V)$ is Galois
invariant.  By \'etale descent it is the base-change of a morphism
$j':D \rightarrow \PP(V)$.  Now $j'$ induces a morphism to
$\widetilde{\PP(V)}$, the blowing up of $\PP(V)$ along $\PP(V_+)\cup
\PP(V_-)$.  Because $j'$ is equivariant for $\iota$ and $\iota_D$, the
quotient morphism $D\rightarrow Z$ factors through $C$, i.e., there is
an induced morphism $i':C\rightarrow Z$.  By a straightforward
enumerative geometry computation, $j'$ has degree $5$ with respect to
$\OO_{\PP(V)}(1)$.  Therefore $i'$ has degree $5$ with respect to
$\OO_Z(1)$.  The degree $3$ multi-section of $\mc{X}_h$ is the image
of $(f,i'):C\rightarrow B\times Z$.

\begin{lem} \label{lem-surj}
% \marpar{lem-surj}   
  If $f$, $g$ and $j$ are general, then $(i')^*:H^0(Z,\OO_Z(1))
  \rightarrow H^0(C,\OO_C(5))$ is surjective.
\end{lem}

\begin{proof}  
  The condition that $(i')^*$ is surjective is an open condition in
  families, hence it suffices to verify $(i')^*$ is surjective for a
  single choice of $f$, $g$ and $j$ -- even one for which
  $\text{Gal}(k(D)/k(B))$ is not $\mathfrak{W}_{3,2}$.  Choose
  homogeneous coordinates $[S_0,S_1]$ on $D$, $[T_0,T_1]$ on $C$ and
  $[U_0,U_1]$ on $B$.  Define $g([S_0,S_1]) = [S_0^2,S_1^2]$ and
  $f([T_0,T_1]) = [T_0^3,T_1^3]$.  Denote by $\mathbf{\mu}_6$ the
  group scheme of $6^\text{th}$ roots of unity.  There is an action of
  $\mathbf{\mu}_6$ on $D$ by $\zeta\cdot[S_0,S_1] = [S_0,\zeta S_1]$.
  This identifies $\mathbf{\mu}_6$ with $\text{Gal}(k(D)/k(B))$.

\medskip\noindent
Let $\mb{e}_{+,0},\mb{e}_{+,1}, \mb{e}_{+,2}$ and
$\mb{e}_{-,0},\mb{e}_{-,1}, \mb{e}_{-,2}$ be ordered bases of $V_+$
and $V_-$ respectively, and let $X_{+,0},X_{+,1},X_{+,2}$ and
$X_{-,0}, X_{-,1}, X_{-,2}$ be the dual ordered bases of $V_+^\vee$
and $V_-^\vee$ respectively.  There is an action of $\mathbf{\mu}_6$
on $V$ by,
$$
\begin{array}{c}
\zeta \cdot [X_{+,0},X_{+,1},X_{+,2},X_{-,0},X_{-,1},X_{-,2}] = \\
\lt[ X_{+,0}, \zeta^2 X_{+,1}, \zeta^4 X_{+,2}, \zeta X_{-,0}, 
\zeta^3 X_{-,1}, \zeta^5 X_{-,2} \rt]
\end{array}
$$ 
and a dual action on $V^\vee$.  Define $j:D\rightarrow \PP(V)$ with
respect to the ordered basis
$\mathbf{e}_{+,0},\dots,\mathbf{e}_{-,2}$, to be the
$\mathbf{\mu}_6$-equivariant morphism,
$$
j([S_0,S_1]) = [S_0^5,S_0^3S_1^2,S_0S_1^3,S_0^4S_1,S_0^2S_1^3,S_1^5].
$$
In this case $U=D_+(U_0U_1)\subset B$ and $\widetilde{W}=W =
D_+(S_0S_1)\subset C$.  It is straightforward to compute $j'$ with
respect to the dual ordered basis $X_{+,0},\dots,X_{-,2}$,
$$
j'([S_0,S_1]) =
[S_1^5,S_0^2S_1^3,S_0^4S_1,S_0S_1^4,S_0S_1^4,S_0^3S_1^2,S_0^5].
$$
As a double-check, observe this is $\mathbf{\mu}_6$-equivariant.  The
induced map $(j')^*$ is,
$$
\begin{array}{ccccccccc}
X_{+,0}X_{+,0} & \mapsto & T_1^5, &
X_{+,0}X_{+,1} & \mapsto & T_0T_1^4, &
X_{+,0}X_{+,2} & \mapsto & T_0^2T_1^3, \\
X_{+,1}X_{+,1} & \mapsto & T_0^2T_1^3, &
X_{+,1}X_{+,2} & \mapsto & T_0^3T_1^2, & 
X_{+,2}X_{+,2} & \mapsto & T_0^4T_1, \\
X_{-,0}X_{-,0} & \mapsto & T_0T_1^4, &
X_{-,0}X_{-,1} & \mapsto & T_0^2T_1^3, &
X_{-,0}X_{-,2} & \mapsto & T_0^3T_1^2, \\
X_{-,1}X_{-,1} & \mapsto & T_0^3T_1^2, &
X_{-,1}X_{-,2} & \mapsto & T_0^4T_1, & 
X_{-,2}X_{-,2} & \mapsto & T_0^5. 
\end{array}
$$
This is surjective by inspection.  
\end{proof}

\subsection{Proof of Theorem~\ref{thm-main}} \label{subsec-proof2}
% \marpar{subsec-proof2}
The subvariety $\mc{X}_h \subset B\times Z$ is a complete intersection
of 3 divisors in the linear system $|\text{pr}_B^*\OO_B(1)\otimes
\text{pr}_Z^* \OO_Z(1)|$, each containing $(f,i')(C)$.  Denote by
$\mc{I}$ the ideal sheaf of $(f,i')(C)\subset B\times Z$, and denote
$I= H^0(B\times Z,\mc{I}\otimes \text{pr}_B^*\OO_B(1)\otimes
\text{pr}_Z^* \OO_Z(1))$.  The projective space of $I$ is the linear
system of divisors on $B\times Z$ in the linear system
$|\text{pr}_B^*\OO_B(1)\otimes \text{pr}_Z^*\OO_Z(1)|$ that contain
$(f,i')(C)$.  The Grassmannian $G'=\text{Grass}(3,I)$ is the parameter
space for deformations of $\mc{X}_h$ that contain $(f,i')(C)$.  For
the same reason as in Corollary~\ref{cor-Sdn}, in $G'$ there is a
countable intersection of dense open subsets parametrizing
subvarieties $\mc{X}'\subset B\times Z$ with $M(K,\mc{X}'_K) \geq 3$
and $I(K,\mc{X}'_K)\mid 4$.  By construction, $\mc{X}'$ contains the
degree $3$ multi-section $(f,i')(C)$.  Therefore $M(K,\mc{X}'_K)=3$
and $I(K,\mc{X}'_K)=1$.  It is straightforward to compute
$\text{pr}_B*[\omega_{\mc{X}'/B}^{\otimes 2}] \cong \OO_B(6)$.  So to
prove the theorem, it suffices to prove every ``very general''
Enriques surface occurs as a fiber of some $\mc{X}'$, i.e., for a
general $[X]\in G$, $X$ occurs as $\text{pr}_Z(\mc{X}'\cap
\pi_B^{-1}(b))$ for some choice of $f, g, i$ and $b\in B$.

\medskip\noindent
A general 0-dimensional, length $3$ subscheme of $Z$ occurs as
$i'(f^{-1}(b))$ for some choice of $f$, $g$, $i$ and $b\in B$.  So for
a general Enriques surface $[X]\in G$ and a general choice of
0-dimensional, length $3$ subscheme of $X$, $X$ is a complete
intersection of $3$ divisors in the linear system $|\OO_Z(1)|$
containing $i'(f^{-1}(b))$ for some choice of $f$, $g$, $i$ and $b$.
To prove that a general $[X]\in G$ is the fiber over $b$ of $\mc{X}'$
for some $f$, $g$, $i$ and $[\mc{X}']\in G'$, it suffices to prove
every divisor in the linear system $|\OO_Z(1)|$ containing
$i'(f^{-1}(b))$ is the fiber over $b$ of a divisor in the linear
system $|\mc{I}\otimes\OO_B(1)\otimes\OO_Z(1)|$.

\medskip\noindent 
There is a short exact sequence,
$$
0 \rightarrow \mc{I}\otimes \text{pr}_Z^*\OO_Z(1) \rightarrow
\text{pr}_Z^*\OO_Z(1)  
\rightarrow \text{pr}_Z^* \OO_Z(1)|_C \rightarrow 0,
$$
giving a short exact sequence,
$$
0 \rightarrow \text{pr}_{B,*}(\mc{I}\otimes \text{pr}_Z^* \OO_Z(1))
\rightarrow \text{pr}_{B,*} \text{pr}_Z^*\OO_Z(1) \rightarrow
\text{pr}_{B,*} (\text{pr}_Z^*\OO_Z(1)|_C) \rightarrow 0.
$$
Because $(i')^*$ is surjective, $\text{pr}_{B,*}(\mc{I}\otimes
\text{pr}_Z^* \OO_Z(1))$ is a locally free sheaf with $h^1=0$.  So it
is $\cong \OO_B^6\oplus \OO_B(-1)^3$.  Twisting by $\OO_B(1)$,
$\text{pr}_{B,*}(\mc{I} \otimes \text{pr}_B^*\OO_B(1)\otimes
\text{pr}_Z^*\OO_Z(1))$ is generated by global sections.  Therefore
every divisor on $Z$ in the linear system $|\OO_Z(1)|$ containing the
scheme $i'(f^{-1}(b))$ is the fiber over $b$ of a divisor on $B\times
Z$ in the linear system $|\mc{I}\otimes \text{pr}_B^*\OO_B(1) \otimes
\text{pr}_Z^*\OO_Z(1)|$.

\bibliography{my}
\bibliographystyle{alpha}

\end{document}